\documentclass[11pt]{amsart}  

 \usepackage[dvips]{epsfig}
\usepackage{graphicx}
\usepackage{latexsym}
\usepackage{amsfonts,amsmath,amssymb}
\newtheorem{corollary}{Corollary}
\newtheorem{lemma}{Lemma}
\newtheorem{example}{Example}
\newtheorem{theorem}{Theorem}
\newtheorem{proposition}{Proposition}
\newtheorem{question}{Question}
\newtheorem{definition}{Definition}
\newtheorem{fact}{Fact}

\begin{document}

\author[M. B\'ona]{Mikl\'os B\'ona}
\title[The absence of a pattern and
the number of occurrences of another]{The absence of a  pattern
 and the number of occurrences of another}
\address{\rm M. B\'ona, Department of Mathematics, 
University of Florida,
358 Little Hall, 
PO Box 118105, 
Gainesville, FL 32611--8105 (USA)
}

\date{}

\begin{abstract} 
Following a question of J. Cooper, we study the expected number
of occurrences of a given permutation pattern $q$ in  permutations that
avoid another given pattern $r$. In some cases, we find the pattern
that occurs least often, (resp. most often) in all $r$-avoiding permutations.
We also prove a few exact enumeration formulae, some of which are
surprising.
\end{abstract}

\maketitle

\section{Introduction}
Let $q=q_1 q_2\ldots q_k$ be a permutation in the
symmetric group $S_k$.  
We say that the permutation
$p=p_1 p_2 \ldots p_n\in
 S_n$ {\it contains a $q$-pattern\/}
 if and only if there is a subsequence
$p_{i_1}p_{i_2}\ldots p_{i_k}$  of $p$ whose elements are in the same
relative order as those in $q$, that is,
$$
p_{i_t}<p_{i_u} \mbox{ if and only if } q_t<q_u
$$
whenever $1\leq t,u\leq k$. If $p$ does not contain $q$, then we say that
$p$ {\em avoids} $q$. 
For example, 41523 contains exactly two occurrences of the pattern 132,
  namely 152 and 153, while 34512 avoids 132. See Chapter 14 of \cite{bona}
for an introduction to pattern avoiding permutations.

It is straightforward to compute, using the linear property of expectation,
 that the average number of $q$-patterns in a randomly selected permutation
of length $n$ is $\frac{1}{k!}{n\choose k}$, where $k$ is the length of $q$.

Joshua Cooper \cite{cooper}
 has raised the following interesting family of questions. 
Let $r$ be a given permutation pattern. What can be said about the
average number of occurrences of $q$ in a randomly selected $r$-avoiding
permutation?

In this paper, we study this family of questions in the case when 
$r=132$. We prove that perhaps surprising result that among patterns of
a fixed length $k$, it is the increasing pattern $12\cdots k$ that occurs
least often and it is the decreasing pattern $k(k-1)\cdots 1$ that occurs
most often in a randomly selected 132-avoiding permutation. 

\section{Preliminaries}
\begin{figure}[ht]
 \begin{center}
  \epsfig{file=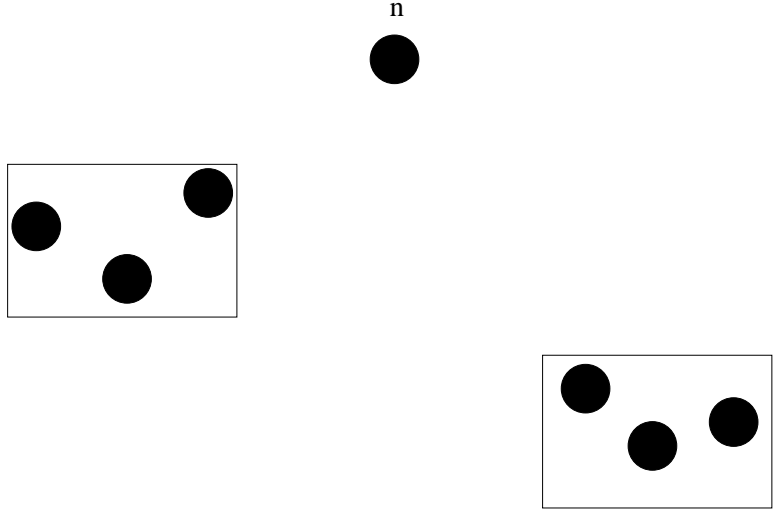}
  \label{generic}
\caption{In a 132-avoiding permutation, all entries preceding the maximal
entry are larger than entries following the maximal entry.}
 \end{center}
\end{figure}

The structure of 132-avoiding permutations is well understood. If 
$p=p_1p_2\cdots p_n$ is such a permutation, and $p_i=n$, then
$p_t>p_u$ must hold for all pairs $(t,u)$ satisfying $t<i<u$. In other
words, all entries on the left of the entry $n$ must be larger than all
entries of the right of $n$. Indeed, if this does not happen, then $p_tnp_u$
is a 132-pattern. This property is so central to the work carried out in this
paper that we illustrate it by Figure \ref{generic}.

Therefore, if $C_n$ denotes the number of 132-avoiding permutations of 
length $n$, then the numbers $C_n$ satisfy the recurrence relation
\begin{equation} 
\label{catalan} C_{n} = \sum_{i=1}^n C_{i-1}C_{n-i},\end{equation}
with $C_0=1$. Hence the numbers $C_n$ are identical to the famous
Catalan numbers $C_n={2n\choose n}/(n+1)$. See Chapter 6 of {\em Enumerative
Combinatorics}\cite{stanley} for a wealth of information on Catalan numbers. 

It follows from the structural property described in the first paragraph
of this section that a 132-avoiding permutation either ends in its largest
entry, or it is {\em decomposable}, that is, it can be cut into two parts
so that every entry that precedes that cut is larger than every entry
that follows that cut. Indeed, if the maximal entry $n$ is not in the
 rightmost position, then one can cut the permutation immediately after $n$
to obtain such a cut. Note that there may be additional ways to cut the 
same permutation.
\begin{example}
The permutation $76834512$ is decomposable. 
One possible cut is  $768|34512$, and another one is  $768345|12$.
\end{example}
This relatively simple structure of 132-avoiding permutations enables 
us to give an exhaustive list of the ways in which a 132-avoiding 
permutation can contain a given pattern. 

\begin{fact} {\em 
 If the 132-avoiding permutation $p$ of length $n$
contains a copy $Q$ of the pattern $q$ of length $k$, then one of 
the following holds.
 (Note
that $q$ itself must be 132-avoiding; otherwise $p$ clearly avoids $q$.) }
\end{fact}
 \begin{enumerate}  
\item If $q$ is not decomposable, that is, if $q$ ends in its largest 
entry, then 
\begin{enumerate}
\item either all of $Q$ must be on the left of $n$,
\item  or all
of $Q$ must be on the right of $n$, 
\item or $Q$ must end in $n$. 
\end{enumerate} 

For instance, if $q=123$, then in $p=678952134$, the subsequence
678 precedes $n=9$, the subsequence 679 ends in 9, and the subsequence
134 follows 9.  

\item If $q$ is decomposable, that is, when $q$ does not end in its largest
entry, and $q$ does not start with its largest entry, then 
\begin{enumerate} 
\item either all of $Q$ is on the left of $n$,
\item or all of $Q$ is on the right of $n$,
\item or the part of $Q$ that precedes a
given cut is on the left of $n$ and the part of $Q$ the follows that cut is
on the right of $n$,
\item  or the part of $Q$ preceding $k$ is on the left of 
$n$, the part of $Q$ following $k$ is on the right of $n$, and the maximal
entry $k$ of $Q$ coincides with $n$.
\end{enumerate}

For instance, if $q=231$, then in $p=786923451$, the subsequence
786 precedes 9, the subsequence 241 follows 9, the subsequence 785 has its
entries 7 and 8 before 9 and its entry 5 after 9 (corresponding to the
cut $23|1$), and the subsequence 895 starts before 9, uses 9, and ends after
9.  

\item  If $q$ is decomposable and $q$ starts with its largest entry $k$,
then 
\begin{enumerate} 
\item either all of $Q$ is on the left of $n$,
\item or all of $Q$ is on the right of $n$,
\item or the part of $Q$ that precedes a
given cut is on the left of $n$ and the part of $Q$ the follows that cut is
on the right of $n$, or 
\item $Q$ starts with $n$, and the rest of $Q$ is on the right of $n$. 
\end{enumerate}

\end{enumerate}
 
\section{Increasing Patterns}

Before proving that among all patterns of length $k$, it is the increasing
pattern $12\cdots k$ that occurs least often in 132-avoiding permutations,
we prove a few general facts about the total number of increasing
patterns in these permutations.

\subsection{A formula for increasing patterns}

Let $a_{n,k}$ be the total number of $12\cdots k$-patterns in all
$C_n$ permutations of length $n$ that avoid 132. 
So for instance, $a_{2,1}=4$, $a_{3,1}=15$, and $a_{2,2}=1$. 
 
Our goal in this subsection is to provide an explicit formula for the
generating function $A_k(x)=\sum_n a_{n,k}x^n$. We will use the well-known
(see for instance Chapter 14 of \cite{bona}) 
explicit formula for the generating function of the Catalan numbers, 
\[C(x)=\sum_{n\geq 0} C_nx^n =  \frac{1-\sqrt{1-4x}}{2x}.\]
We will prove the following theorem.

\begin{theorem} \label{increasing} We have
\begin{equation}
\label{forone} 
A_1(x)=\sum_{n\geq 1}nC_nx^n =\sum_{n\geq 1}{2n\choose n}x^n  - \sum_{n\geq 1}
C_nx^n =\frac{1}{\sqrt{1-4x}} -\frac{1-\sqrt{1-4x}}{2x}.\end{equation}
 Furthermore, 
for all positive integers $k\geq 2$, we have
\begin{equation}
\label{explicit} A_k(x)= A_1(x) \left(\frac{xC(x)}{1-2xC(x)}\right)^{k-1}=
A_1(x)\left( \frac{1}{2\sqrt{1-4x}} - \frac{1}{2}\right)^{k-1}=
A_1(x)F^{k-1}(x).
\end{equation}
\end{theorem}

\begin{proof}
For $k=1$, the claim is obvious, since an increasing subsequence of length
one is just an entry of a permutation. 

For larger $k$, an increasing subsequence of length $k$ is an indecomposable
pattern. Hence the ways in which it can occur in the 132-avoiding
permutation $p$
are listed in case (1) of Fact 1.  This leads
to the recurrence relations \[a_{n,k}=2\sum_{i=1}^n a_{i-1,k}C_{n-i}
+\sum_{i=1}^n a_{i-1,k-1}C_{n-i},\]
or in terms of generating functions, 
\begin{eqnarray} \label{monrec} A_k(x) & = & 2xA_k(x)C(x)+xA_{k-1}(x)C(x) \\
& = & A_{k-1}(x) \frac{xC(x)}{1-2xC(x)},
\end{eqnarray}
and our claim follows by induction on $k$. 
\end{proof}

Note that in particular, (\ref{explicit}) implies that for $1\leq k< l$
we have
\begin{equation} \label{stationary}
A_k(x)A_l(x) = A_{k+1}(x)A_{l-1}(x).
\end{equation}

\subsection{Why the Increasing Pattern is Minimal}
For a given pattern $q$, let $t_n(q)$ denote the number of all 
occurrences of the pattern $q$ in all 132-avoiding permutations of 
length $n$. So in particular, if $q$ is the increasing pattern $12\cdots k$,
then $t_n(q)=a_{n,k}$. 

The main result of this section is the following theorem, which shows that
 no pattern of a given length occurs less often in the set of
132-avoiding permutations than the increasing pattern of that length. 

\begin{theorem} \label{minimal}
Let $q$ be any pattern of length $k$. Then for all
positive integers $n$, we have
$t_q(n)\geq a_{n,k}$.
\end{theorem}

Before proving the theorem, we need to introduce some simple machinery
to simplify notation.

\begin{definition}
Let $G(x)=\sum_{n\geq 0}g_nx^n$ and $H(x)=\sum_{n\geq 0}h_nx^n$
be two power series. We say that $G(x)\leq H(x)$ if $g_n\leq h_n$ for
all $n\geq 0$.
\end{definition}

\begin{proposition} \label{multiply}
Let $G(x)$, $H(x)$ and $W(x)$ be three power series with non-negative
real coefficients so that $G(x)\leq H(x)$ holds. Then 
\[G(x)W(x) \leq H(x)W(x).\]
\end{proposition}

\begin{proof} 
The coefficient of $x^n$ in $H(x)W(x)-G(x)W(x)$ is 
\[\sum_{i=0}^n (h_i-g_i)w_{n-i},\] which is a sum of non-negative real
numbers, and is hence non-negative.
\end{proof}

We can now return to the proof of Theorem \ref{minimal}.

\begin{proof} (of Theorem \ref{minimal})
We prove the statement by induction on $k$. For $k=1$, the statement is
obvious.

Now let us assume that the statement is true for all positive integers
less than $k$, and prove it for $k$. We distinguish three cases. Each of
these cases will be handled by analyzing recurrence relations, which 
may sometimes seem somewhat cumbersome. Therefore, at the beginning of
each case, we will give an intuitive description of that case. 

In a permutation $p=p_1p_2\cdots p_n$, we say that $i$ is a {\em descent}
if $p_i>p_{i+1}$. Otherwise, we say that $i$ is an {\em ascent}. 

An overview of the cases is as follows. First, we treat patterns ending in
their largest entry. Then we treat patterns that contain only one descent,
say in position $j$. Finally, we treat all remaining patterns, comparing 
the patterns whose first descent is in position $j$ to the pattern whose
{\em only} descent is in position $j$.

See Figure \ref{threesteps} for an illustration.

\begin{figure}[ht]
 \begin{center}
  \epsfig{file=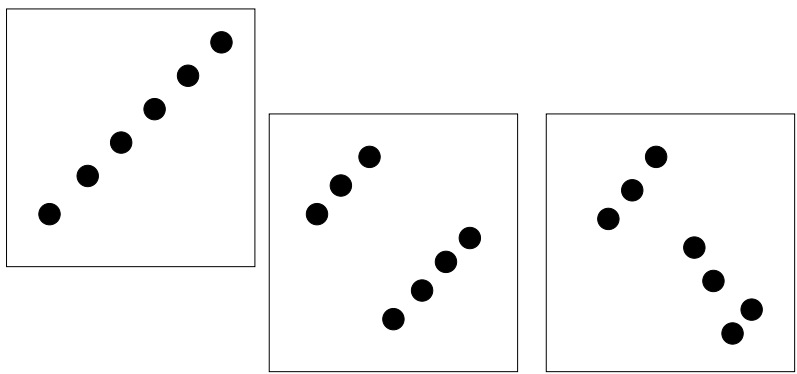}
  \label{threesteps}
\caption{Three types of patterns that we compare.}
 \end{center}
\end{figure}

\begin{enumerate}
\item[Case 1] When $q$ ends in its largest entry $k$. 
Let $q'$ be the pattern
obtained from $q$ by removing the largest entry $k$ from the end of $q$.
We then show that $q$ is more frequent than $12\cdots k$ by 
showing that $q$ can be contained in a 132-avoiding permutation in the
same ways as $12\cdots k$, and then using the induction hypothesis.
In other words, the difference between $a_{n,k}$ and $t_n(q)$ is caused
by whatever happens before the last position of these patterns.

The possible ways in which $q$ can occur in a 132-avoiding permutation $p$
are listed in case (1) of Fact 1. 
Therefore, we have the recurrence relation
\[t_n(q)=2\sum_{i=1}^n t_{i-1}(q)C_{n-i} + \sum_{i=1}^n t_{i-1}(q')C_{n-i},\]
leading to the generating function identities
\[T_q(x)=2xC(x)T_q(x) + xC(x)T_{q'}(x),\]
\begin{equation} \label{compare}
T_q(x)=T_{q'}(x)\frac{xC(x)}{1-2xC(x)}.\end{equation}
Comparing formulae (\ref{monrec}) and (\ref{compare}), we see that
$T_q(x)$ is obtained from   $T_{q'}(x)$ by the same operation as
$A_k(x)$ is obtained from $A_{k-1}(x)$, namely by a multiplication by
the power series $F(x)=\frac{xC(x)}{1-2xC(x)}$. 
As \[[x^n]T_{q'}(x)=t_n(q') \geq a_{n,k-1}=[x^n]A_{k-1}(x)\] by our
induction hypothesis, and $F(x)=\sum_{n\geq 1}{2n-1\choose n-1}x^n$ has
non-negative coefficients, our claim is immediate
 by Proposition \ref{multiply}.
\item[Case 2] Let us now consider the case in which $q$ does not 
end in its largest entry $k$, and the only descent of $q$ is 
in the position in which $k$ occurs. We will subsequently
 see that all remaining
cases will easily reduce to this one. Let us
say that $k$ is in the $j$th position in $q$, with $j<k$.
 That is, 
$q=q_{k,j}=(k-j+1)(k-j+2)\cdots k123 (k-j)$. For instance, $q_{7,3}
=5671234$. We will show that 
\begin{equation} \label{firstcompare}
t_n(q_{k,j}) \geq t_n(12\cdots k)=a_{n,k}.
\end{equation}

The main idea is the following. The pattern $q_{k,j}$ looks very similar
to the increasing pattern, hence the ways in which  $q_{k,j}$ can be contained
in a 132-avoiding permutation are also similar to the ways in which the
increasing pattern can. So the numbers $t_n(q_{k,j})$ and $a_{n,k}$ satisfy
very similar recurrence relations, and where they differ, they differ in
the way we predicted.

\begin{enumerate}
\item[Subcase 2.1] Let us assume first that $j\neq 1$, that is, that $k$ is
not in the first position in $q_{k,j}$.
As $q_{k,j}$ is decomposable by a cut after the $j$th position, the ways in
which $q_{k,j}$ can be contained in a 132-avoiding permutation are described
in case (2) of Fact 1. That list leads to the recurrence relation
\begin{eqnarray} 
\label{recfort} t_n(q_{k,j}) & = & 2\sum_{i=1}^n C_{i-1}t_{n-i}(q_{k,j})
 + \sum_{i=1}^n 
a_{i-1,j}a_{n-i,k-j}\\
& + & \sum_{i=1}^n 
a_{i-1,j-1}a_{n-i,k-j},  \end{eqnarray}
and hence
\[T_{q_{k,j}}(x)=2xC(x)T_{q_{k,j}}(x) + xA_{j}(x)A_{k-j}(x)
 +  xA_{j-1}(x)A_{k-j}(x), \]
\begin{equation}
\label{expltq} 
T_{q_{k,j}}(x)=\frac{xA_{k-j}(x)(A_{j}(x)+A_{j-1}(x))}{1-2xC(x)}.
\end{equation}
We need to show that $T_{q_{k,j}}(x)\geq A_k(x)$. Comparing formulae
(\ref{explicit}) and (\ref{expltq}), that is equivalent to 
\[\frac{xA_{k-j}(x)(A_{j}(x)+A_{j-1}(x))}{1-2xC(x)}\geq A_1(x)F(x)^{k-1},\]
or \begin{equation} \label{secondcase}
\frac{x}{1-2xC(x)} A_1(x)^2 (F(x)^{k-3} + F(x)^{k-2})
 \geq A_1(x)F(x)^{k-1}.\end{equation}
Inequality (\ref{secondcase}) will be proved if we can show that
 \begin{equation} 
\label{stronger} \frac{x}{1-2xC(x)} A_1(x)(1+F(x) ) \geq F(x)^2.
\end{equation}
Indeed, (\ref{stronger}) implies (\ref{secondcase}) by Proposition 
\ref{multiply}, choosing $W(x)=F^{k-3}(x)$. 
On the other hand,  (\ref{stronger}) is equivalent to 
\[\frac{x}{2(1-4x)^{3/2}} + \frac{x}{2(1-4x)} -\frac{1}{4(1-4x)} +
\frac{1}{4} \geq \frac{1}{4(1-4x)} -\frac{1}{2\sqrt{1-4x}} +
\frac{1}{4} ,\]
\[\frac{x}{2(1-4x)^{3/2}} +\frac{x-1}{2(1-4x)}
 + \frac{1}{2\sqrt{1-4x}} \geq 0.\]
The coefficient of $x^n$ on the left-hand side is 0 if $n=0$, $n=1$, or
$n=2$,  
and is $b_n={2n-3\choose n-2}(2n-1)-3\cdot 2^{2n-3}  +{2n-1\choose n-1}$
if $n\geq 3$. If we replace $n$ by $n+1$, the negative summand in the
above expression of  $b_n$, that is,
$3\cdot 2^{2n-3}$, grows fourfold, whereas a routine computation shows that 
the sum of the two positive terms grows $\frac{4n^2+14n+6}{n^2+3n+2}$-fold.
This fraction is larger than 4 for all $n\geq 3$, showing that $b_n\geq 0$
for all $n$, and  our claim is proved.

\item[Subcase 2.2]
If $j=1$, then a minor modification is necessary since if a copy of
$q$ contains $n$, then it has to start with $n$. Hence formula
 (\ref{recfort}) becomes 
\begin{eqnarray} 
\label{recforj1} t_n(q_{k,1}) & = & 2\sum_{i=1}^n C_{i-1}t_{n-i}(q_{k,1})
 + \sum_{i=1}^n 
a_{i-1,1}a_{n-i,k-1}\\
& + & \sum_{i=1}^n 
C_{i-1}a_{n-i,k-1},  \end{eqnarray}
 So only the last sum is
different
from what it was in (\ref{recfort}).
This leads to the generating function identities
\[T_{q_{k,1}}(x)=2xC(x)T_{q_{k,1}}(x) + 
xA_{k-1}(x)(A_1(x)+C(x)), \]
\begin{equation}
\label{explfor1} T_{q_{k,1}}(x)=\frac{xA_{k-1}(x)(A_1(x)+C(x))}{1-2xC(x)}.
\end{equation}
Comparing formulae (\ref{monrec}) and (\ref{explfor1}), the inequality
$A_k(x) \leq  T_{q_{k,1}}(x)$ is now proved by Proposition \ref{multiply},
since  $C(x) \leq A_{1}(x)+C(x)$.
\end{enumerate}

\item[Case 3] Finally, there is the case when $k$ is in the $j$th position of
$q$ for some $j<k$, but $j$ is not the only descent of $q$.  
We claim that then copies of $q$ occur even more frequently than copies
of $q_{k,j}$, roughly because even in segments where $q_{k,j}$ is increasing,
$q$ is not. 

That is, we will prove that
\begin{equation}
\label{ncompare}
t_n(q_{k,j}) \leq t_n(q).
\end{equation}
This, together with (\ref{firstcompare}) will complete the proof of 
Theorem \ref{minimal}. 

As $q$ is decomposable by a cut after its $j$th position, the ways in which
$q$ can be contained in a 132-avoiding permutation are described in case
(2) of Fact 1. Let $q^{<1>}$ denote the pattern formed by the first $j$ 
entries of $q$, and let $q^{<2>}$ be the pattern formed by the remaining
$k-j$ entries of $q$. Then we have the recurrence relation
\[t_n(q) \geq  2\sum_{i=1}^n C_{i-1}t_{n-i}(q) + \sum_{i=1}^n 
t_{i-1}(q^{<1>}) t_{n-i} (q^{<2>}) + \sum_{i=1}^n 
t_{i-1}(q^{<1'>}) t_{n-i} (q^{<2>}).  \]
Here $q^{<1'>}$ is the pattern obtained from $q^{<1>}$ by removing its last
(and also largest) entry. 
Note that $t_n(q)$ is at least as large as the right-hand side, and not 
necessarily equal to it. That is because, unlike $q_{k,j}$, the pattern $q$
may be decomposable by other cuts, in addition to the cut after its $j$th
entry. The existence of such cuts would add extra summands to the
right-hand side. 

The last displayed inequality leads to the generating function inequalities
\[T_q(x) \geq 2T_q(x)xC(x) + xT_{q^{<1>}}(x)T_{q^{<2>}}(x) + 
xT_{q^{<1'>}}(x)T_{q^{<2>}}(x),\] 
\[T_q(x)(1-2xC(x)) \geq xT_{q^{1}}(x)T_{q^{2}}(x) + 
xT_{q^{<1'>}}(x)T_{q^{<2>}}(x).\]
Note that as $1/(1-2xC(x))=1/\sqrt{1-4x}=\sum_{n\geq 0}{2n\choose n}x^n$
has non-negative coefficients, the last displayed inequality remains true
if we multiply both sides by $1/(1-2xC(x))$. This leads to the inequality
\begin{equation} \label{explgen}
T_q(x) \geq \frac{ xT_{q^{<2>}}(x)(T_{q^{<1>}} + 
T_{q^{<1'>}}(x))}{1-2xC(x)}.\end{equation}
We can now compare formulae (\ref{expltq}) and (\ref{explgen}). We see that
by our induction hypothesis,
 each factor on the right-hand side of (\ref{explgen}) is at least as
large as the corresponding factor on the right-hand side of (\ref{expltq}).
That is, $xT_{q^{<2>}}(x)\geq xA_{k-j}(x)$, since $q^{<2>}$ is a pattern of
length $k-j$, and $T_{q^{<1>}} + 
T_{q^{<1'>}}(x)\geq A_{k}(x) + A_{k-1}(x)$ by a summand-wise comparison. 
Hence, by Proposition \ref{multiply}, we have that
$T_q(x) \geq T_{q_{k,j}}(x)$, and
the proof of Theorem \ref{minimal} is complete.
\end{enumerate}
 \end{proof}

\section{Why the decreasing pattern is maximal}
In this section we prove that the decreasing pattern $k(k-1)\cdots 21$
occurs more frequently in 132-avoiding permutations than any other pattern
of length $k$. The structure of the proof will be very similar to 
that of the minimality of the increasing pattern, but there will be 
more technical difficulties. 

Let $d_{n,k}$ denote the number of decreasing subsequences of length $n$
in all 132-avoiding permutations of length $n$. Then we have 
\[d_{n,1}=a_{n,1}=nC_n=\frac{n}{n+1}{2n\choose n}.\]
For larger values of $k$, consider the set of all $C_n$ permutations
of length $n$ that avoid 132. In that set, for every $1\leq j\leq k-1$,
 there are
\[\sum_{i=1}^n d_{i-1,j}d_{n-i,k-j}\] copies of $k(k-1)\cdots 1$ in which
the first $j$ entries are on the left of $n$, and the last $k-j$ entries
are on the right of $n$. (The index $i$ denotes the position of the entry 
$n$ in a permutation of length $n$.)
 In addition, there are $\sum_{i=1}^n C_{i-1}d_{n,k-1}$ copies of
$k(k-1)\cdots 1$ that start with the entry $n$. Finally, there are the
$2\sum_{i=1}^n C_{i-1}d_{n-i,k}$ 
copies of $k(k-1)\cdots 1$ that are either entirely on the left of $n$, 
or entirely on the right of $n$.
This leads to the recurrence relation
\begin{equation} \label{decrec}
d_{n,k}=\sum_{j=1}^{k-1}
\sum_{i=1}^n d_{i-1,j}d_{n-i,k-j}+\sum_{i=1}^n C_{i-1}d_{n,k-1}+
 2\sum_{i=1}^n C_{i-1}d_{n-i,k}
\end{equation}
and the generating function identities
\[
D_k(x) = 2xC(x)D_k(x) + xC(x)D_{k-1}(x) +\sum_{j=1}^{k-1} xD_j(x)D_{k-j}(x), 
\] 
\begin{equation}
\label{recd} D_k(x) = \frac{ xC(x)D_{k-1}(x) +
\sum_{j=1}^{k-1} xD_j(x)D_{k-j}(x)}{1-2xC(x)}.
\end{equation}

The following Corollary provides an estimate for the ``growth'' of the
power series $D_k(x)$. It is worth comparing its result
 with Theorem \ref{increasing}.

\begin{corollary} \label{easycor}
We have \begin{equation}
\label{for2} D_2(x)=\frac{xD_1(x)}{(1-4x)}, \end{equation}
and 
\begin{equation}
\label{larger} 
D_k(x)\geq \frac{xD_{k-1}(x)}{(1-4x)} 
\end{equation}
for $k\geq 3$.
\end{corollary}

\begin{proof}
The first displayed identity (the special case of $k=2$)
immediately follows from (\ref{recd}) if we recall that $1-2xC(x)=
\sqrt{1-4x}$ and that $D_1(x)+C(x)=\sum_{n\geq 0}{2n\choose n}x^n
=\frac{1}{\sqrt{1-4x}}$.

The general formula (\ref{larger}) follows from (\ref{recd}) if
we remove all summands from the right-hand side except for
$xC(x)D_{k-1}(x)$ and $xD_1(x)D_{k-1}$, (we can do this since all the
removed terms have non-negative coefficients), and then again, we
recall that $D_1(x)+C(x)=\frac{1}{\sqrt{1-4x}}$.
\end{proof}

The following Lemma is a natural counterpart of its much simpler analogue
(\ref{stationary}). It shows that the sequence of power series $D_1(x),D_2(x),
\cdots $ is log-convex in a certain sense. 

\begin{lemma} \label{logc}
For all positive integers $2\leq a \leq b$ we have
\[D_a(x)D_b(x) \leq D_{a-1}(x)D_{b+1}(x).\]
\end{lemma}

\begin{proof}
Induction on $a+b$. The smallest value of $a+b$ for which the statement is
not trivial is 4. The non-trivial statement then is that
$D_2(x)^2\leq D_1(x)D_3(x)$. By (\ref{easycor}), this is equivalent
to $D_1(x)^2\frac{x^2}{(1-4x)^2} \leq D_1(x)D_3(x)$. In order to prove
the latter, it suffices to show that 
\[D_1(x)\frac{x^2}{(1-4x)^2}\leq D_3(x),\] and that is immediate by
(\ref{for2}).

Now let us assume that the statement holds for $a+b=m-1$ and prove it 
for $a+b=m$. By (\ref{recd}), it suffices to show that 
\[D_a(x)C(x)D_{b-1}(x)+\sum_{j=1}^{b-1}D_a(x)D_j(x)D_{b-j}(x) 
\leq \] \[\leq 
 D_{a-1}C(x)D_b(x) + \sum_{j=1}^{b}D_{a-1}(x)D_j(x)D_{b+1-j}(x).
\]
Note that the right-hand side has one more summand than the left-hand side.
Now we carry out a series of pairwise comparisons. First, for the two terms
preceding the summation signs, we have that $D_a(x)C(x)D_{b-1}(x)\leq
D_{a-1}(x)C(x)D_b(x)$ since $C(x)$ has non-negative coefficients, and
by our induction hypothesis, $D_{a}(x)D_{b-1}(x)\leq D_{a-1}(x)D_{b}(x)$.
For the terms after the summation signs, for $j<a$, we have that
\[D_{a}(x)D_j(x)D_{b-j}(x) \leq D_{a-1}(x)D_j(x)D_{b+1-j}(x)\] since our
induction hypothesis implies that $D_{a}(x)D_{b-j}(x)\leq
D_{a-1}(x)D_{b+1-j}(x)$. The induction hypothesis applies since 
$a+b-j<a+b$.

For $a\leq j \leq b-1$, we claim that 
\[D_{a}(x)D_j(x)D_{b-j}(x) \leq D_{a-1}(x)D_{j+1}(x)D_{b-j}(x).\]
(That is, we skip one summand, and we compare the $j$th summand of the 
left-hand side to the $(j+1)$st summand to the right-hand side.)
Indeed, our induction hypothesis implies that 
$D_{a}(x)D_j(x)\leq D_{a-1}(x)D_{j+1}(x)$. The induction hypothesis applies 
since $a+j<a+b$.

With this, each summand of the left-hand side was injectively associated to
a weakly larger summand of the right-hand side, proving our claim.  
\end{proof}

Now we are in a position to state and prove the main result of this section.

\begin{theorem} \label{dectheo}
Let $q$ be a pattern of length $k$. Then the inequality
\[t_n(q)\leq d_{n,k} \]
holds.
\end{theorem}

\begin{proof}

We prove the statement by induction on $k$ and $n$. We know that the statement
holds for $k=1$ and it is routine to verify it for $k=2$.
 Let us now assume that it is true for all 
patterns shorter than $k$. Let us further assume that for patterns of 
length $k$, the statement holds for all permutations shorter than $n$.
(The initial cases of $n<k$ are obvious.) Let us now prove that the
statement hold for permutations of length $n$, and patterns of length $k$. 
 
Again, our proof proceeds by cases. We first handle patterns that start with
their largest entry, then patterns which contain only one ascent, in position
$j$. Finally, we cover the remaining cases, comparing patterns whose first
ascent is in position $j$ to the pattern whose only ascent is in position $j$.
See Figure \ref{threedown} for an illustration of some of these cases.

\begin{figure}[ht]
 \begin{center}
  \epsfig{file=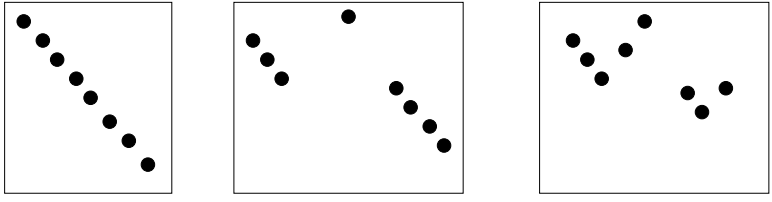}
  \label{threedown}
\caption{Three kinds of patterns we compare.}
 \end{center}
\end{figure}

\begin{enumerate}
\item[Case 1] Let us first consider the case when
 $q$ starts with its largest entry $k$. In this case, simply use the
fact that the pattern obtained from $q$ by removing its first entry occurs 
less often than the pattern $k(k-1)\cdots 21$, and it is decomposable by
no more cuts than  $k(k-1)\cdots 21$.

 In this case, $q$ is
decomposable since it can be cut right after its first entry. Let us
say that $q$ is decomposable with $v$ distinct cuts, and let 
$q^{1,f}$ and $q^{1,b}$ denote the patterns before and after the first cut,
let $q^{2,f}$ and $q^{2,b}$ denote the patterns before and after the second
cut, and so on, up to $q^{v,f}$ and $q^{v,b}$ for the patterns on the two
sides of the last cut. (The letters $f$ and $b$ stand for ``front'' and
 ``back''.)

Note that the total length of  $q^{j,f}$ and $q^{j,b}$ is always $k$ for
every $j$, and that  $|q^{1,f}|=1$  and $|q^{1,b}|=k-1$.

Then we have the recurrence relation
\begin{equation}
\label{easycase} t_n(q)=
\sum_{j=1}^v \sum_{i=1}^n t_{i-1}(q^{j,f})t_{n-i}(q^{j,b})
+ \sum_{i=1}^n C_{i-1}t_{n-i}(q^{1,b})+ 2\sum_{i=1}^nC_{i-1}t_{n-i}(q) 
\end{equation}

It is now straightforward to compare $d_{n,k}$ and $t_n(q)$ by comparing
the corresponding summands of recurrence relations (\ref{decrec}) and 
(\ref{easycase}). Let us first compare the two double sums. In those
sums, $j$ indexes the cuts of the respective patterns. As the decreasing
pattern of length $k$ has $k-1$ cuts, $j$ ranges from 1 to $k-1$ in
 (\ref{decrec}). In (\ref{easycase}), $j$ ranges from 1 to $v$, where
$v\leq k-1$ is the number of cuts that $q$ has. So the first double sum has
more terms. We also claim that even the terms that the second double sum
does have are smaller than the corresponding terms in the first double sum.

Indeed, if $|q^{j,f}|=y$ and $|q^{j,b}|=k-y$, then $t_{i-1}(q^{j,f})\leq
d_{i-1,y}$ and $t_{n-i}(q^{j,b})\leq d_{n-i,k-y}$ by our induction hypothesis,
so $t_{i-1}(q^{j,f})t_{n-i}(q^{j,b})\leq d_{i-1,y} d_{n-i,k-y}$.

Comparing the second sums of  (\ref{decrec}) and 
(\ref{easycase}) is even simpler. Their summands agree in the term $C_{i-1}$,
and by our induction hypothesis, we know that $t_{n-i}(q^{1,b})\leq 
d_{n-i,k-1}$. Finally, comparing the third sums of (\ref{decrec}) and 
(\ref{easycase}) we use the fact that $n-i<n$, so by our induction hypothesis,
$t_{n-i}(q)<d_{n-i,k}$. 

\item[Case 2]
 When $q$ has only one ascent and $q$ does not start with its largest
entry $k$. Just as in the previous section,
this is the heart of the proof. The remaining cases will easily reduce to
this one. 
Let $q^{k,h}$ be the pattern of length $k$ that avoids 132 and has
only one ascent, in position $h$. As $q$ does not start in its largest
entry $k$, this means that $k$ must be the $(h+1)$st entry of $q$. That is, 
\[q^{k,h}=(k-1)\cdots (k-h)k(k-h-1)\cdots 321.\]
For instance $q^{5,2}=43521$.

The difficulty of this case is that there is a way in which $q^{k,h}$
can be contained in a 132-avoiding permutation in which $k(k-1)\cdots 21$
cannot, namely by having $h$ entries on the left of $n$ and $k-h-1$ entries
on the right of $n$. As we will see, the ways in which the decreasing pattern
can be contained in such a permutation and $q^{k,h}$ cannot are more 
prevalent, but that is not obvious. We will need Lemma \ref{logc} to prove
that fact. 

The structure of $q^{k,h}$ is somewhat similar to that of the decreasing 
pattern. That is,  $q^{k,h}$ can be cut after each entry starting
with $k$ (in position $h+1$). Cutting immediately after position $j$ (with 
$j\geq h+1$)
will result in the two patterns $q^{j,h}$ and the decreasing pattern
of length $k-j$. This leads to the recurrence relation
 
\begin{equation} \label{decrec2}
t_n(q^{k,h})=\sum_{j=h+1}^{k-1} \sum_{i=1}^n t_{i-1}(q^{j,h})d_{n-i,k-j}
+ \sum_{i=1}^n d_{i-1,h}d_{n-i,k-h-1} +  
2\sum_{i=1}^n C_{i-1}t_{n-i}(q^{k,h}).
\end{equation}

We can carry out a pairwise comparison of corresponding sums in  formulae
 (\ref{decrec}) and (\ref{decrec2}). This is easiest for the third sums
in those two formulae: indeed, as $n-i<n$, our induction hypothesis
implies that $t_{n-i}(q^{k,h})\leq d_{n-i,k}$ for all $k$, $i$, and
$h$, so the third sum appearing in    (\ref{decrec}) is larger than the
third sum appearing in  (\ref{decrec2}).

As far as the double sums are concerned, for $j=h+1,h+2,\cdots ,k-1$,
our induction hypothesis implies that $t_{i-1}(q^{j,h})\leq
 d_{i-1,j}$. Therefore, each term of the double sum of  (\ref{decrec2})
is at most as large as the corresponding term of (\ref{decrec}).

Therefore, the claim $t_n(q^{k,h})\leq d_{n,k}$ will be proved if we can
show that the remaining sum in (\ref{decrec2}) is less than the remaining
sums in (\ref{decrec}),  that is, that  
\[\sum_{i=1}^n d_{i-1,h}d_{n-i,k-h-1} \leq 
\sum_{j=1}^{h}
\sum_{i=1}^n d_{j,i-1}d_{n-i,k-j}+\sum_{i=1}^n C_{i-1}d_{n-i,k-1}.\]

We show the stronger statement that the above inequality remains true
even if we remove all summands from the double sum in which $j\neq 1$.
(Note that as $k$ is not in the first position of $q$, we know that $h-1\geq
1$, so this will leave a non-empty set of summands.) 
In other words, we claim that 
  
\[\sum_{i=1}^n d_{i-1,h}d_{n-i,k-h-1} \leq 
\sum_{i=1}^n d_{i-1,1}d_{n-i,k-1}+\sum_{i=1}^n C_{i-1}d_{n,k-1}.\]
This is equivalent to the generating function inequalities
\[xD_h(x)D_{k-h-1}(x)\leq xD_{1}(x)D_{k-1}(x) + xC(x)D_{k-1}(x).\]
\[D_h(x)D_{k-h-1}(x) \leq D_{k-1}(x) (1-4x)^{-1/2}.\]
By Lemma \ref{logc}, we know that $D_h(x)D_{k-h-1}(x)\leq D_1(x)D_{k-2}(x)$,
so it suffices to prove that
\[ D_1(x)D_{k-2}(x) \leq D_{k-1}(x) (1-4x)^{-1/2}.\]
By Corollary \ref{easycor}, we know that
$\frac{x}{1-4x} D_{k-2}(x) \leq D_{k-1}(x)$, so our claim will be proved
if we can show that \[D_1(x)\leq \frac{x}{(1-4x)^{3/2}}.\]
The last displayed inequality holds since 
$[x^n]D_1(x)=\frac{n}{n+1}{2n\choose n}$ whereas
 $[x^n]\frac{x}{(1-4x)^{3/2}}=(2n-1){2n-2\choose n-1}$. A routine computation
shows that the latter is larger as soon as $2<n+1$, or $1<n$.

\item Finally, we consider the case when $q$ has more than one ascent, and
$q$ does not start with its maximal entry.
Let us assume that the maximal entry $k$ of
 $q$ is in the $h$th position. We will
show that then $t_n(q)\leq t_n(q^{k,h})$. The main idea behind the proof is
that $q^{k,h}$ is decomposable at every place where $q$ is, and after 
decomposition, its parts are close to the decreasing pattern.

Let us inductively assume that we know the statement $t_n(q)\leq t_n(q^{k,h})$
for all patterns shorter than $k$.

The set of positions after  which $q$ can be cut is a subset of the
set of positions after which $ q^{k,h}$ can be cut. Just as in Case (1),
let us  say that $q$ is decomposable with $v$ distinct cuts, and let 
$q^{1,f}$ and $q^{1,b}$ denote the patterns before and after the first cut,
let $q^{2,f}$ and $q^{2,b}$ denote the patterns before and after the second
cut, and so on, up to $q^{v,f}$ and $q^{v,b}$ for the patterns on the two
sides of the last cut. Note that $|q^{1,f}|=h$    and $|q^{1,b}|=k-h$.

Then, similarly to Case (1), we  have the recurrence relation
\begin{equation} \label{decrec3}
t_n(q)=\sum_{j=1}^v \sum_{i=1}^n t_{i-1}(q^{j,f})t_{n-i}(q^{j,b})
+ \sum_{i=1}^n t_{i-1}(q^{1,f'})t_{n-i}(q^{1,b})+
 2\sum_{i=1}^nC_{i-1}t_{n-i}(q),
\end{equation}
where $q^{1,f'}$ is the pattern $q^{1,f}$ with its last (and largest) entry
removed. 

The inequality $t_n(q)\leq t_n(q^{k,h})$ is now obvious by pairwise comparing
the three summands in (\ref{decrec2}) and (\ref{decrec3}), and using the
induction hypotheses. (In particular, when comparing the first summands, 
we use the fact that if $|q^{j,h}|=|q^{j,f}|$, then $t_{i-1}(q^{j,h})\geq
t_{i-1}(q^{j,h})$ by the induction hypothesis made in this case. Recall
that we assumed that for patterns $q$ shorter than $k$, it is true that
$t_n(q)\leq t_n(q^{|q|,h})$, where $h$ is the first ascent of $q$.)
\end{enumerate}
\end{proof}

\section{Asymptotic Enumeration}
Theorem \ref{increasing} provides an explicit formula for $A_2(x)$.
From that formula, it is routine to deduce that
\begin{equation} \label{asc2ex}
a_{n,2}=2^{2n-1} - \frac{1}{4}{2n+2\choose n+1} - \frac{n}{2n+1}{2n\choose n}
\end{equation}
for $n\geq 1$.

Similarly, an explicit formula for $D_2(x)$ can be obtained by setting
$k=2$ in (\ref{recd}). That is, 
\begin{equation}
D_2(x)=\frac{xC(x)D_1(x)+xD_1(x)D_1(x)}{1-2xC(x)}=\frac{x}{(1-4x)^{3/2}}
-\frac{1}{2(1-4x)}+ \frac{1}{2\sqrt{1-4x}}.
\end{equation}
This yields
\begin{equation}
\label{dec2ex} 
d_{n,2}={2n-1\choose n-1}(n+1) -2^{2n-1} 
\end{equation}
for $n=1$.

Comparing formulae (\ref{asc2ex}) and (\ref{dec2ex}) and using 
Stirling's formula, we see that in 
132-avoiding permutatations of length $n$, there are $c\sqrt{n}$ times
as many inversions as non-inversions.

For larger values of $k$, it would be perhaps tempting to think 
that increasing patterns of length $k$ become very rare, since they are
the rarest of all patterns of length $k$. However, Theorem \ref{increasing}
shows that $a_{n,k}$ is of exponential order 4, since the dominant singularity
of $A_k(x)$ is at $1/4$. In fact, Theorem \ref{increasing} implies
that  $A_k(x)$ is a sum of various powers of $(1-4x)$, the dominant of which 
is $(1-4x)^ {-\frac{k}{2} }$.
It then follows from singularity analysis that $a_{n,k}\sim c \cdot 4^n 
n^{\frac{k-2}{2}}$.

\section{Further Directions}

A simple analysis of the proofs of Theorems \ref{minimal} and \ref{dectheo}
show that the inequalities $a(n,k) \leq t_n(q) \leq d(n,k)$ are sharp is
$n$ is large enough compared to $k$. 

We have seen that if we consider 132-avoiding permutations, then among all
patterns of length $k$, the increasing pattern is the least likely to
occur and the decreasing pattern is the most likely to occur. This suggests
the following natural question.

\begin{question}
Let $r$ be any pattern, and let $t_{r,q}(n)$ be the number of all copies
of $q$ in all $r$-avoiding permutations of length $n$. Let us assume that
among all patterns $q$ of length $k$, it is the increasing pattern that
minimizes $t_{r,q}(n)$ for all $n$. 

Is it then true that among all patterns $q$ of length $k$, it is the 
decreasing pattern that maximizes $t_{r,q}(n)$?
\end{question}

Another direction of research is the following. 

\begin{question}
Let $q_1$ and $q_2$ be two patterns of the same length, and assume that
for some positive integer $N$, the inequality
\[t_{r,q_1}(N) < t_{r,q_2}(N) \]
holds. Is it then true that 
\[t_{r,q_1}(n) < t_{r,q_2}(n) \] for all $n>N$?
\end{question}

In other words, is it true that the relation between the frequency of
$q_1$ and $q_2$ in $r$-avoiding permutations depends only on $r$, $q_1$ and
$q_2$, or does it depend on $n$ as well? 

Lemma \ref{logc} and formula \ref{stationary} show interesting combinatorial
properties of the power series $D_1(x), D_2(x),\cdots $ and  
$A_1(x), A_2(x),\cdots $. These properties are easy to express in terms
of combinatorial objects, without power series. However, is there a 
combinatorial proof for them?

\end{document}